\documentclass[a4paper,12pt]{article}
\usepackage{amsmath,amsthm,amssymb}

\theoremstyle{plain}
\newtheorem*{theorem*}{Theorem}
\newtheorem{theorem}{Theorem}[section]
\newtheorem{lemma}[theorem]{Lemma}
\newtheorem{corollary}[theorem]{Corollary}
\newtheorem{proposition}{Proposition}[subsection]
\theoremstyle{remark}
\newtheorem{example}[proposition]{Example}

\numberwithin{equation}{subsection}

\DeclareMathOperator{\centr}{center}

\begin{document}

\title{Non-rational divisors over non-Gorenstein terminal 
singularities}
\author{D.~A. Stepanov\thanks{The work was partially supported by
RFBR grant no. 02-01-00441, and Grants of Leading Scientific Schools
no. 489.2003.1. and no. 1910.2003.1}}
\date{}

\maketitle

\begin{abstract}
Let $(X,o)$ be a germ of a 3-dimensional terminal singularity of
index $m\geqslant 2$. If $(X,o)$ has type $cAx/4$, $cD/3-3$, 
$cD/2-2$, or $cE/2$, then assume that the standard equation of
$X$ in $\mathbb{C}^4/\mathbb{Z}_m$ is non-degenerate with respect
to its Newton diagram. Let $\pi\colon Y\to X$ be a resolution. We
show that there are not more than $2$ non-rational divisors $E_i$,
$i=1,2$, on $Y$ such that $\pi(E_i)=o$ and discrepancy 
$a(E_i,X)\leqslant 1$. When such divisors exist, we describe them
as exceptional divisors of certain blowups of $(X,o)$ and study
their birational type.
\end{abstract}

\section{Introduction}
In this paper, we continue the study of resolutions of terminal
singularities started in \cite{cD} and \cite{DEndeg}.

Let $(X,o)$ be a germ of a 3-dimensional terminal singularity
defined over the field $\mathbb{C}$ of complex numbers. Consider
a resolution $\pi\colon Y\to X$ and let $E\subset Y$ be a prime 
divisor such that $\pi(E)=o$ and discrepancy $a(E,X)\leqslant 1$. 
Note that if $\pi$ is a divisorial resolution, then $E$ does
exist (see \cite{Kaw}, \cite{Mar}). On the other hand, the number of
such divisors is finite (here we identify two divisors over $X$ if
they give the same discrete valuations of the field $k(X)$). 

What can be said about the birational type of the algebraic surface
$E$? It is known that $E$ is birationally ruled (\cite{C3f}, 
Corollary~2.14). Moreover, if the singularity $(X,o)$ is of type 
$cA/m$, $m\geqslant 1$, then the surface $E$ is rational 
(\cite{E1blowups}, Proposition~2.4). When the singularity $(X,o)$ is
of type $cD$, the surface $E$ is either rational or birationally 
isomorphic to $\mathbb{P}^1\times C$, where $C$ is a (hyper)elliptic
curve. If this non-rational divisor $E$ exists, it is unique
(\cite{cD}). When $(X,o)$ is a general singularity of type $cE$, 
the non-rational divisor with low discrepancy is again unique
and birational to the surface $\mathbb{P}^1\times C$, but the
curve $C$ can be non-hyperelliptic (\cite{DEndeg}).

In this paper, we study the case when $(X,o)$ is a general
non-Gorenstein (i.~e., the canonical divisor $K_X$ is not a Cartier
divisor) 3-dimensional terminal singularity. By ``general'' we
mean the following. Any non-Gorenstein terminal singularity is
analitically isomorhic to one of singularities of 
Theorem~\ref{T:classification} which we call standard. The 
singularity $(X,o)$ is general if its standard equation in 
$\mathbb{C}^4/\mathbb{Z}_m$ is non-degenerate with respect to its 
Newton diagram. 
\begin{theorem*}
Let $\pi\colon Y\to X$ be a resolution of 3-dimensional 
non-Gorenstein terminal singularity $(X,o)$. If $(X,o)$ is of type
$cAx/4$, $cD/3-3$, $cD/2-2$, or $cE/2$, then additionally assume
that the standard defining equation of $X$ is non-degenerate
with respect to its Neqton diagram. Then there are not more than
$2$ non-rational divisors $E_i$, $i=1,2$, such that $\pi(E_i)=o$
and discrepancy $a(E_i,X)\leqslant 1$.
\end{theorem*}
In all the cases when the non-rational divisors exist, we describe
them as exceptional divisors of certain blowups of the singularity
$(X,o)$ and study their birational type. 

In section~\ref{S:prelim}, we recall the analylic classification
of 3-dimensinal non-Gorenstein terminal singularities and state
some lemmas useful for working with discrepancies and resolutions. 
In section~\ref{S:proof}, we prove our Theorem by case-by-case 
analysis of all types of non-Gorenstein terminal singularities. We 
do not consider the case of $cA/m$-singularities because it was 
completely studied by Yu.~G. Prokhorov in~\cite{E1blowups}.

\section{Preliminaries}\label{S:prelim}
Let the cyclic group $\mathbb{Z}_m$ act on the space $\mathbb{C}^n$ 
as follows: $x_i\to \varepsilon^{a_ir}x_i$, $i=1,\dots,n$, where
$x_i$ are the coordinates in $\mathbb{C}^n$, $\varepsilon$ is a
primitive $m$-th root of unity, $a_i\in\mathbb{Z}$, 
and $r\in\mathbb{Z}_m$ is a residue modulo $m$. We shall denote
the quotient space $\mathbb{C}^n/\mathbb{Z}_m$ by $\mathbb{C}^n/
\mathbb{Z}_m(a_1,a_2,\dots,a_n)$ or by 
$\frac{1}{m}(a_1,a_2,\dots,a_n)$.

The classification of 3-dimensional non-Gorenstein terminal
singularities was obtained by Danilov, Mori, Koll\'{a}r, and 
Shephard-Barron.
\begin{theorem}\label{T:classification}
(\cite{Mo}) Let $X$ be a germ of a 3-dimensional terminal 
singularity of index $\geqslant 2$. Then there is an embedding
of $X$ to $\mathbb{C}^4/\mathbb{Z}_m$ such that one of the 
following holds: 
\begin{description}
\item[(cA/m)] 
$X\simeq\{xy+f(z,u)=0\}\subset\frac{1}{m}(\alpha,-\alpha,1,0)$ where
$\alpha$ is an integer prime to $m$ and $f(z,u)\in\mathbb{C}\{z,u\}$
is a $\mathbb{Z}_m$-invariant.
\item[(cAx/4)]
$X\simeq\{x^2+y^2+f(z,u)=0\}\subset\frac{1}{4}(1,3,1,2)$ where
$f(z,u)\in\mathbb{C}\{z,u\}$ is a $\mathbb{Z}_4$-semi-invariant
and $u\notin f(z,u)$ (the coeficient of the monomial $u$ in the
power series $f$ is zero).
\item[(cAx/2)]
$X\simeq\{x^2+y^2+f(z,u)=0\}\subset\frac{1}{2}(0,1,1,1)$ where
$f(z,u)\in(z,u)^4\mathbb{C}\{z,u\}$ is a $\mathbb{Z}_2$-invariant.
\item[(cD/3)]
$X\simeq\{\varphi(x,y,z,u)=0\}\subset\frac{1}{3}(1,2,2,0)$ where
$\varphi$ has one of the following forms:
  \begin{description}
  \item[(cD/3-1)]
  $\varphi=u^2+x^3+yz(y+z)$,
  \item[(cD/3-2)]
  $\varphi=u^2+x^3+yz^2+xy^4\lambda(y^3)+y^6\mu(y^3)$ where 
$\lambda(y^3)$, $\mu(y^3)\in\mathbb{C}\{y^3\}$ and 
$4\lambda^3+27\mu^2\ne 0$,
  \item[(cD/3-3)]
  $\varphi=u^2+x^3+y^3+xyz^3\alpha(z^3)+xz^4\beta(z^3)+
yz^5\gamma(z^3)+z^6\delta(z^3)$ where $\alpha(z^3)$, $\beta(z^3)$,
$\gamma(z^3)$, $\delta(z^3)\in\mathbb{C}\{z^3\}$.
  \end{description}
\item[(cD/2)]
$X\simeq\{\varphi(x,y,z,u)=0\}\subset\frac{1}{2}(1,1,0,1)$ where 
$\varphi$ has one of the following forms:
  \begin{description}
  \item[(cD/2-1)]
  $\varphi=u^2+xyz+x^{2a}+y^{2b}+z^c$ where $a$, $b\geqslant 2$,
$c\geqslant 3$,
  \item[(cD/2-2)]
  $\varphi=u^2+y^2z+\lambda yx^{2a+1}+g(x,z)$ where 
 $\lambda\in\mathbb{C}$, $a\geqslant 1$, $g(x,z)\in(x^4,x^2z^2,z^3)
\mathbb{C}\{x,z\}$.
  \end{description}
\item[(cE/2)]
$X\simeq\{u^2+x^3+g(y,z)x+h(y,z)=0\}\subset\frac{1}{2}(0,1,1,1)$ where
$g(y,z)\in(y,z)^4\mathbb{C}\{y,z\}$, $h(y,z)\in
(y,z)^4\mathbb{C}\{y,z\}\setminus(y,z)^5\mathbb{C}\{y,z\}$.
\end{description}
The index of $X$ is equal to the order of the cyclic group 
$\mathbb{Z}_m$.
\end{theorem}
\begin{theorem}
(\cite{KoShB}) Let $X$ be one of the hyperquotient singularities
$$\{\varphi(x,y,z,u)=0\}\subset\mathbb{C}^4/\mathbb{Z}_m$$
listed in Theorem~\ref{T:classification}. Assume that 
$\varphi(x,y,z,u)=0$ defines an isolated singularity at $0$ and the 
action of $\mathbb{Z}_m$ is free on $X$ outside $0$. Then $X$ is 
terminal.
\end{theorem}

Let $f=f(x_1,x_2,\dots,x_n)$ be a covergent power series such that
$f(0)=0$ and $(\{f=0\},0)\subset$$(\mathbb{C}^n,0)$ is an isolated
singularity. We denote by $\Gamma(f)$ the Newton diagram of the 
series $f$. If $f$ is non-degenerate with respect to its Newton
diagram (in the sequel, we say simply that $f$ is non-degenerate), 
then there is a Varchenko-Hovanski\u{\i} embedded toric resolution
of the singularity $(\{f=0\},0)$ (see \cite{Varchenko}).
Moreover, if the group $\mathbb{Z}_m$ acts on $\mathbb{C}^n$ and
$f$ is its semi-invariant, then we can repeat the construction from
\cite{Varchenko} and obtain an embedded toric resolution of the
quotient singularity 
$$(X,o)=(\{f=0\},0)/\mathbb{Z}_m\subset\mathbb{C}^n/\mathbb{Z}_m\,.$$
Here all necessary toric varieties and morphisms are built with 
respect to the lattice $N'$ dual to the lattice $M'$ of monomials 
invariant under the action of $\mathbb{Z}_m$, $M'\subset\mathbb{Z}^n$. 
This easy observation was pointed out to us by S.~A. Kudryavtsev.

Recall that the embedded toric resolution $\pi\colon Y\to X$ of the
singularity $(X,o)$ is determined by a certain subdivision of
the non-negative octant $\mathbb{R}_{\geqslant 0}^{n}$.  If $\Sigma$ is 
the corresponding fan, then let $\widetilde{\mathbb{C}^n}=$
$X(\Sigma,N')$ be the toric variety built from $\Sigma$ and let
$\tilde\pi\colon\widetilde{\mathbb{C}^n}\to
\mathbb{C}^n/\mathbb{Z}_m$ be the natural birational morphism. 
Then $\pi$ is the restriction of the morphism $\tilde\pi$ to the
proper transform $Y$ of the singularity $X$.

Exceptional divisors of the morphism $\tilde\pi$ are in one-to-one
correspondence with 1-dimensional cones of the fan $\Sigma$. Take
a 1-dimensional cone $\tau$, its exceptional divisor 
$E_{\tau}\subset\widetilde{\mathbb{C}^n}$, and let 
$E_{\tau}|_Y=\sum m_jE_j$. Further, let $w=(w_1,\dots,w_n)$ be the 
primitive vector of the lattice $N'$ along the cone $\tau$. 
The diagram $\Gamma(f)$ lies in the space $(\mathbb{R}^n)^*$ dual
to $\mathbb{R}^n=\mathbb{R}\otimes\mathbb{Z}^n$; we denote the
corresponding pairing by $\langle \cdot ,\cdot \rangle$. Now we want 
to calculate discrepancy $a(E_j,X)$.
\begin{lemma}\label{L:discrepancy}
$a(E_j,X)=m_j(w_1+w_2+\dots+w_n-1-w(f))$, where $w(f)=$
$=\min\{\langle w,v\rangle\,|\,v\in\Gamma(f)\}$.
\end{lemma}
\begin{proof}
Arguing as in \cite{Varchenko}, \S10, we find an affine neighborhood
$U\simeq\mathbb{C}^n$ of the generic point $E_\tau$ in 
$\widetilde{\mathbb{C}^n}$ with
coordinates $y_1,\dots,y_n$ such that the equation $y_1=0$ 
defines $E_\tau\cap U$ and the morphism
$\tilde\pi|_U\colon U\to\mathbb{C}^n/\mathbb{Z}_m$ is given by
the formulae: 
$$x_1=y_{1}^{w_1}y_{2}^{a^{2}_{1}}\dots y_{n}^{a^{n}_{1}}\,,$$
$$\dots\dots$$
$$x_n=y_{1}^{w_n}y_{2}^{a^{2}_{n}}\dots y_{n}^{a^{n}_{n}}$$
for some $a^i=(a^{i}_{1},\dots,a^{i}_{n})\in N'\cap
\mathbb{R}^{n}_{\geq 0}$. To prove the lemma, it remains only
to lift the differential form $dx_1\wedge\dots\wedge dx_n$ to $U$
and to apply the adjunction formula.
\end{proof}
\begin{corollary}\label{C:discr}
If $a(E_j,X)\leqslant 1$, then $w_1+\dots+w_n-1-w(f)\leqslant 1$.
\end{corollary}

Note that if the vectors $w$, $e_1$, $e_2$,$\dots$,$e_n$, where
$e_i=(0,\dots,\underset{i}{1},\dots,0)$, generate the lattice $N'$,
then the exceptional divisors $E_j$ are birationally isomorphic to 
the divisors $E_{w,j}$ respectively, $\sum E_{w,j}=E_w|_{X_w}$, 
where $X_w$ is the proper transform of $X$ under the weighted blowup
$$\nu_w\colon\mathbb{C}_{w}^{n}\to\mathbb{C}^n/\mathbb{Z}_m\,.$$ 
This follows from the fact that for any two subdivisions of the 
non-negative octant $\mathbb{R}_{\geqslant 0}^{n}$ there is a common
subsubdivision. The exceptional divisor $E_w$ of $\nu_w$ is 
isomorphic to the weighted projective space 
$\mathbb{P}(w_1,\dots,w_n)$.
The divisor $\sum E_{w,j}$ is defined in $\mathbb{P}(w_1,\dots,w_n)$
by the equation
$$f_{\rho(w)}(x_1,\dots,x_n)=0\,,$$
where $f_{\rho}$ corresponds to the face 
$$\rho(w)=\{v\in\Gamma(f)\,|\,\langle w,v\rangle=w(f)\}\,,$$
$$f_{\rho(w)}=
\sum_{(m_1,\dots,m_n)\in\rho(w)}
a_{m_1\dots m_n}x_{1}^{m_1}\dots x_{n}^{m_n}$$
$$\text{ if }
f=\sum_{\substack{(m_1,\dots,m_n)\in\\ \mathbb{Z}_{\geqslant 0}}}
a_{m_1\dots m_n}x^{m_1}\dots x_{n}^{m_n}\,.$$
Now the letters $x_1,\dots,x_n$ denote the
quasihomogeneous coordinates in the space 
$\mathbb{P}(w_1,\dots,w_n)$. We often use this abuse of notation in 
the sequel; the wright meaning of the letters is clear from the 
context.

Now suppose that the vectors $w$, $e_1$,$\dots$,$e_n$ generate
some sublattice $N''\subset N'$. Consider the subdivision of
the octant $\mathbb{R}_{\geqslant 0}^{n}$ by the vector $w$, i.~e.,
the fan $\Sigma_w$ consisting of the cones $\sigma_i=
\langle e_1,\dots,\underset{i}{w},\dots,e_n\rangle$ and all their
faces. So obtained morphism 
$$\mu_w\colon\widetilde{\mathbb{C}^{n}_{w}}\to
\mathbb{C}^n/\mathbb{Z}_m$$ 
is not a weighted blowup. We shall call it a \emph{pseudo blowup
with the weight} $w$. It is easily proved that its 
exceptional divisor $\Tilde E_w\simeq\mathbb{P}(w_1,\dots,w_n)/G$, 
where $G=N'/N''$ is a cyclic group, and the equation of
$\sum \Tilde E_{w,j}$ is the same as above.

Let $(X,o)$ be one of the terminal singularities listed in 
Theorem~\ref{T:classification}, let $\nu_w$ be its weighted
blowup or pseudo blowup, and let $E_w$ be the exceptional divisor
of the morphism $\nu_w\colon X_w\to X$. Denote by $E'$ the surface
in $\mathbb{P}(w_1,\dots,w_n)$ covering $E_w$ (if $\nu_w$ is a
weighted blowup, then $E'=E_w$).
\begin{lemma}\label{L:rationality}
Suppose that the surface $E'$ is irreducible and has only rational
singularities. Then the surface $E_w$ is rational.
\end{lemma}
\begin{proof}
We can consider $E'$ as a divisor over some terminal $cDV$-point.
Take a resolution $\pi\colon\widetilde E'\to E'$ of singularities of
the surface $E'$. According to \cite{C3f}, Corollary~2.14, $E'$ 
is birationally ruled. Thus
$P_2(\widetilde E')=h^0(2K_{\widetilde E'})=0$. On the other hand,
$E'$ is a hypersurface in the space $\mathbb{P}(w_1,\dots,w_n)$,
hence $h^1(\mathcal{O}_{E'})=0$. Since $E'$ has only rational
singularities, we have $h^1(\mathcal{O}_{\widetilde E'})=
h^1(\mathcal{O}_{E'})=0$. Therefore $\widetilde E'$ is rational by
Castelnuovo's criterion and thus $E_w$ is rational too.
\end{proof}

If the blowup $\nu$ has a non-rational exceptional divisor, we shall
sometimes say that the blowup $\nu$ is \emph{non-rational}.

\section{Proof of Theorem}\label{S:proof}
\subsection{Terminal singularities of type $cAx/4$}\label{S:cAx4}
Consider the singularity $(X,o)$ of type $cAx/4$, i.~e.,
$$X\simeq\{\varphi=x^2+y^2+f(z,u)=0\}\subset\frac{1}{4}(1,3,1,2)\,,$$ 
where $f(z,u)\in\mathbb{C}\{z,u\}$ is a $\mathbb{Z}_4$-semi-invariant
and $u\notin f(z,u)$. We assume that the defining series $\varphi$ 
is non-degenerate. Then the singularity $(X,o)$ has an embedded
toric resolution $\pi\colon Y\to X$. Divisors with center at $o$ and
discrepancy $a\leqslant 1$ belong to any divisorial resolution of
$X$, so if there is a non-rational divisor $E$ over $(X,o)$, 
$\centr(E)=o$, $a(E,X)\leqslant 1$, then $E$ belongs to the 
resolution $\pi$. We saw in section~\ref{S:prelim} that $E$ is
birationally isomorphic to the exceptional divisor (or to its
irreducible component) of some weighted blowup or pseudo blowup
$\nu_w$. Thus we can suppose that $E$ is given in 
$\mathbb{P}(w_1,w_2,w_3,w_4)$ (or in 
$\mathbb{P}(w_1,w_2,w_3,w_4)/G$) by the part $\varphi_w$ of the 
series $\varphi$. It is clear that if $E$ is non-rational, then
the polynomial $\varphi_w$ contains at least one of the monomials 
$x^2$ or $y^2$. But if it contains both of them, i.~e., 
$\varphi_w=x^2+y^2+f_w(z,u)$, then in the affine chart $u\ne 0$ 
the surface $E_w$ is defined as
$$\{x^2+y^2+f_w(z,1)=0\}\subset\mathbb{C}^4/G_1\,,$$
and in the chart $z\ne 0$ as 
$$\{x^2+y^2+f_w(1,u)=0\}\subset\mathbb{C}^4/G_2\,,$$
where $G_1$, $G_2$ are finite cyclic subgroups of 
$\mathbf{GL}_\mathbb{C}(3)$.
Now it is obvious that $E_w$ has only rational singularities,
thus by Lemma~\ref{L:rationality} $E_w$ is rational.
So, we can assume that $\varphi_w=x^2+f_w(z,u)$ or 
$\varphi_w=y^2+f_w(z,u)$.

Since $\varphi$ is a $\mathbb{Z}_4(1,3,1,2)$-semi-invariant and the
singularity $(X,o)$ is isolated,
the series $f$ contains the monomial $u^{2n+1}$. If $n$ is the 
minimal number with the property $u^{2n+1}\in f$, we say that
the singularity $(X,o)$ is of type $cA_{2n}x/4$. Now let us find
all blowups and pseudo blowups $\nu$ of the singularity $(X,o)$ of 
type $cA_{2n}x/4$ such that the exceptional divisor of $\nu$ can
be non-rational and has discrepancy $a\leqslant 1$.

We have to find all primitive vectors $w\in\mathbb{Z}^4+
\frac{1}{4}(1,3,1,2)\mathbb{Z}$ such that \\
either (i) $w_1+w_2+w_3+w_4-1-2w_1\leqslant 1$, $w_2>w_1$,
$(2n+1)w_4\geqslant 2w_1$, \\
or (ii) $w_1+w_2+w_3+w_4-1-2w_2\leqslant 1$, $w_1>w_2$,
$(2n+1)w_4\geqslant 2w_2$.
\begin{proposition}\label{P:cAx4}
Let the primitive vector $w\in\mathbb{Z}^4+\frac{1}{4}(1,3,1,2)
\mathbb{Z}$ satisfy conditions (i) or (ii). Then $w$ is one of the
following: \\
1) $\frac{1}{4}(4k+1,4k+3,1,2)$, $k\leqslant n/2$, $k\in 
\mathbb{Z}_{\geqslant 0}$; \\
2) $\frac{1}{4}(4k+3,4k+5,3,2)$, $k\leqslant (n-1)/2$, $k\in
\mathbb{Z}_{\geqslant 0}$; \\
3) $\frac{1}{4}(4k+5,4k+3,1,2)$, $k\leqslant (n-1)/2$, $k\in
\mathbb{Z}_{\geqslant 0}$; \\
4) $\frac{1}{4}(4k+3,4k+1,3,2)$, $k\leqslant n/2$, $k\in
\mathbb{Z}_{\geqslant 0}$.
\end{proposition}
\begin{proof} It is an easy arithmetic calculation. For example,
assume (i). Since $w_2>w_1$, the inequality for discrepancy has
the form $w_3+w_4<2$. Taking into account that 
$w_3\in\frac{1}{4}\mathbb{Z}$, $w_4\in\frac{1}{2}\mathbb{Z}$, and
$2w_3\equiv w_4 \mod \mathbb{Z}$,
we get the following possibilities: \\
$w_4=1/2,\quad w_3=1/4,\,3/4,\,5/4\,;$ \\
$w_4=1,\quad w_3=1/2\,;$ \\
$w_4=3/2,\quad w_3=1/4\,.$ \\
Assume $(w_3,w_4)=(1/4,1/2)$. Then we have 
$w_2-w_1+3/4-1\leqslant 1$,
i.~e., $w_2-w_1\leqslant 5/4$. On the other hand, $w_1\leqslant
(2n+1)w_4=n/2+1/4$. If we combine these inequalities with
$w_1\equiv w_3\mod\mathbb{Z}$, we obtain $w_1=\frac{1}{4}(4k+1)$,
$w_2=\frac{1}{4}(4k+3)$, $k\leqslant n/2$, i.~e., case 1).

Let us also consider the possibility $(w_3,w_4)=(5/4,1/2)$. It
follows that $w_2-w_1\leqslant 1/4$. But this is impossible because
the difference $w_2-w_1$ is always multiple of $1/2$.

Other cases can be done in a similar way. In the sequel, we omit
such calculations. 
\end{proof}

Note that vectors 1)--4) give weighted blowups but not
pseudo blowups. 

The exceptional divisor $E_1$ of the blowup 
$\nu_1=\frac{1}{4}(4k+1,4k+3,1,2)$
(weighted blowup with the weight $\frac{1}{4}(4k+1,4k+3,1,2)$) 
is defined by the equation 
$$\{x^2+f_{2k+\frac{1}{2}}(z,u)=0\}\subset
\mathbb{P}(4k+1,4k+3,1,2)\,.$$
If $E_1$ is non-rational, it is irreducible and reduced.
Then discrepancy 
$$a(E_1,X)=(1/4)(4k+1+4k+3+1+2)-1-2k-1/2=1/4\,.$$
It is obvious that $E_1$ is a cone over the hyperelliptic curve 
$C=\{x^2+f_{2k+\frac{1}{2}}(z,u)=0\}\subset\mathbb{P}(4k+1,1,2)$. 
Genus of a curve in a weighted projective plane can be found by
methods from \cite{Dolgachev}. Genus of the curve $C$ is 
$g(C)\leqslant 2k$.

The exceptional divisor $E_2$ of the blowup $\nu_2=
\frac{1}{4}(4k+3,4k+5,3,2)$ is defined by the equation
$$\{x^2+f_{2k+\frac{1}{2}}(z,u)=0\}\subset
\mathbb{P}(4k+3,4k+5,3,2)\,.$$
If it is irreducible and reduced, its discrepancy $a(E_2)=3/4$.
The divisor $E_2$ is a cone over a hyperelliptic curve of genus 
$$g\leqslant
\begin{cases}
   2m-1,\quad k=3m\,, \\
   2m+1,\quad k=3m+1\,, \\
   2m+2,\quad k=3m+2\,.
\end{cases}
$$

The exceptional divisor $E_3$ of the blowup $\nu_3=
\frac{1}{4}(4k+5,4k+3,1,2)$ is defined by the equation
$$\{y^2+f_{2k+\frac{1}{2}}(z,u)=0\}\subset
\mathbb{P}(4k+5,4k+3,1,2)\,.$$
If $E_3$ is irreducible and reduced, discrepancy $a(E_3)=1/4$.
The surface $E_3$ is a cone over a hyperelliptic curve of genus 
$g\leqslant 2k+1$.

The exceptional divisor $E_4$ of the blowup $\nu_4=
\frac{1}{4}(4k+3,4k+1,3,2)$
is defined by the equation 
$$\{y^2+f_{2k+\frac{1}{2}}(z,u)=0\}\subset
\mathbb{P}(4k+3,4k+1,3,2)\,.$$
If $E_4$ is irreducible and reduced, discrepancy $a(E_4)=3/4$. 
Divisor $E_4$ is a cone over a hyperelliptic curve of genus
$$
g\leqslant
\begin{cases}
   2m,\quad k=3m\,,\\
   2m+1,\quad k=3m+1\text{ or } k=3m+2\,.
\end{cases}
$$

It is clear that the blowups $\nu_1$ and $\nu_3$, $\nu_2$ and 
$\nu_4$ can not simultaneously be non-rational. Indeed, assume
for example that $\nu_1$ is non-rational. It follows
that for the weights $w(z)=1$, $w(u)=2$ the function
$f$ has the weight $w(f)=8k_1+2$. But if $\nu_3$ is also 
non-rational, then $w(f)=8k_2+6$, contradiction.
Other pairs of blowups can be non-rational. 
\begin{example}
Consider the singularity
$$\{x^2+y^2+z^{18}+z^6u^6+u^{15}=0\}\subset\frac{1}{4}(1,3,1,2)$$
of type $cA_{14}x/4$. Make the blowups $\nu_1=\frac{1}{4}(9,11,1,2)$
and $\nu_2=\frac{1}{4}(15,17,3,2)$. The exceptional divisor 
$$E_1\colon\{x^2+z^{18}+z^6u^6=0\}\subset\mathbb{P}(9,11,1,2)$$
of the first one is a cone over a singular curve of genus $g=2$. 
The exceptional divisor 
$$E_2\colon\{x^2+z^6u^6+u^{15}=0\}\subset\mathbb{P}(15,17,3,2)$$
of the second blowup is a cone over a singular curve of genus $g=1$.
\end{example}

We see that there is not more than $2$ non-rational divisors
with discrepancy $a\leqslant 1$ over a non-degenerate singularity
of type $cAx/4$.

\subsection{Terminal singularities of type $cAx/2$}\label{S:cAx2}
Consider the singularity $(X,o)$ of type $cAx/2$, i.~e.,
\begin{equation}\label{E:cAx2}
X\simeq\{x^2+y^2+f(z,u)=0\}\subset\frac{1}{2}(0,1,1,1)\,,
\end{equation}
where $f(z,u)\in(z,u)^4\mathbb{C}\{z,u\}$ is a 
$\mathbb{Z}_2$-semi-invariant. Here our proof does not depend on
the fact whether given singularity is non-degenerate. Following
\cite{Hay}, \S8, assume that if the weights of variables are
$w(z)=w(u)=1/2$, then the weight $w(f)$ of the series $f$ equals $k$.
If $k$ is even, we make the weighted blowup $\nu_0=
\frac{1}{2}(k,k+1,1,1)$. If $k$ is odd, we make the blowup
$\nu_1=\frac{1}{2}(k+1,k,1,1)$. We shall only consider $\nu_0$, 
the other case can be done in a similar way. 

We have $\nu_0\colon\widetilde{\mathbb{C}^4}\to\mathbb{C}^4/
\mathbb{Z}_2(0,1,1,1)$ and the variety $\widetilde{\mathbb{C}^4}$
is coversd by $4$ affine charts. In the first one $U_1\simeq
\frac{1}{k}(1,-1,-1,-1)$, the proper transform $\Tilde X$ of the
singularity $X$ is given by the equation
$$1+xy^2+f_k(z,u)+x(\dots)=0\,.$$
It is clear that in $U_1$ the variety $\Tilde X$ is non-singular.

In the second chart $U_2\simeq\frac{1}{k+1}(1,1,-1,-1)$,
$$\Tilde X\cap U_2\colon x^2+y+f_k(z,u)+y(\dots)=0\,.$$
Here $\Tilde X$ is non-Gorenstein at the origin only, where it has
a cyclic terminal quotient singularity of type 
$\frac{1}{k+1}(1,-1,-1)$. The third and the fourth charts are
isomorphic to $\mathbb{C}^4$. In the third one
the variety $\Tilde X\cap U_3$ is defined by the equation
$$x^2+y^2z+f_k(1,u)+z(\dots)=0\,.$$
Since $(X,o)$ is an isolated singularity, singularities of
$\Tilde X\cap U_3$ lie only on the exceptional divisor $\{z=0\}$. 
It is obvious that all of them are isolated $cDV$-points. 
Similarly, in the fourth chart the variety $\Tilde X$ has only
isolated $cDV$-points.

Let $E$ be the exceptional divisor of the bloup $\nu_0$ of the
singularity $(X,o)$. We have
$$E\simeq\{x^2+f_k(z,u)=0\}\subset\mathbb{P}(k,k+1,1,1)\,.$$
If it is non-rational, it is irreducible and reduced, 
discrepancy $a(E,X)=(1/2)(k+k+1+1+1)-1-k=1/2$, and the surface
$E$ is a cone over a hyperelliptic curve of genus $g\leqslant k-1$.

Take an arbitrary resolution $\pi\colon Y\to\Tilde X\overset{\nu_0}
{\to} X$. All non-rational divisors with discrepancy $a\leqslant 1$
appear in $\pi$. But $cDV$-singularities of the variety 
$\Tilde X$ produce only divisors with discrepancies $a(E_i,\Tilde X)
\geqslant 1$, so that $a(E_i,X)>1$. Any resolution of the cyclic
quotient singularity from the second chart of $\Tilde X$ contains
with discrepancies $\leqslant 1$ only rational divisors. Thus
$E$ is the unique non-rational divisor with $a\leqslant 1$ over
the singularity $(X,o)$. We have proved the following
\begin{proposition}\label{P:cAx2}
Any resolution of the singularity $(X,o)$ of type $cAx/2$ contains
not more than $1$ non-rational divisor $E$ with discrepancy 
$a(E,X)\leqslant 1$ and $\centr_X(E)=o$. Let $(X,o)$ be defined
by the equation~\eqref{E:cAx2}.
Then the non-rational divisor $E$ can be realized as the exceptional 
divisor of the weighted blowup $\nu_0=\frac{1}{2}(k,k+1,1,1)$ 
(if $k$ is even), or as the exceptional divisor of the weighted 
blowup $\nu_1=\frac{1}{2}(k+1,k,1,1)$ (if $k$ is odd).
In both cases $E$ is a cone over a hyperelliptic curve of genus 
$g\leqslant k-1$.
\end{proposition}

\begin{example}
Consider the singularity 
$$\{x^2+y^2+z^6+u^6=0\}\subset\frac{1}{2}(0,1,1,1)$$
and its weighted blowup $\frac{1}{2}(4,3,1,1)$. Its exceptional
divisor
$$E\simeq\{y^2+z^6+u^6=0\}\subset\mathbb{P}(4,3,1,1)$$
is a cone over a curve of genus $2$.
\end{example}

\subsection{Terminal singularities of type $cD/3$}
\subsubsection{$cD/3-1$}
Consider the singularity $(X,o)$ of type $cD/3-1$, i.~e.,
$$X\simeq\{u^2+x^3+yz(y+z)=0\}\subset\frac{1}{3}(1,2,2,0)\,.$$
This singularity can be resolved by an explicit calculation.
There are no non-rational divisors with discrepancy $a\leqslant 1$ 
over $(X,o)$.

\subsubsection{$cD/3-2$}\label{S:cD32}
Consider the singularity $(X,o)$ of type $cD/3-2$, i.~e., 
\begin{equation}\label{E:cD32}
X\simeq\{u^2+x^3+yz^2+xy^4\lambda(y^3)+y^6\mu(y^3)=0\}
\subset\frac{1}{3}(1,2,2,0)\,,
\end{equation}
where $\lambda(y^3)$, $\mu(y^3)\in\mathbb{C}\{y^3\}$ and 
$4\lambda^3+27\mu^2\ne 0$. Note that the last condition guarantees
that the singularity $(X,o)$ is non-degenerate. However, we shall not 
use this fact. We shall proceed as in case $cAx/2$ in 
section~\ref{S:cAx2}.

Consider the weighted blowup $\nu=\frac{1}{3}(2,1,4,3)$ (see
\cite{Hay}, \S9) of the given singularity. It can be easily
verified that in the first, in the second, and in the fourth charts
the blouwn up variety $\Tilde X$ is non-singular. In the third
chart $U_3\simeq\frac{1}{4}(2,3,3,1)$
$$\Tilde{X_3}=\Tilde X\cap U_3\simeq
\{u^2+x^3+yz+\lambda_0xy^4+\mu_0y^6+z(\dots)=0\}\,.$$
At the origin the variety $\Tilde{X_3}$ has a singularity 
analytically isomorphic to 
$$\{u^2+y^2+z^2+x^3=0\}\subset\frac{1}{4}(2,3,3,1)\,.$$
It is odvious that it has type $cAx/4$ and is non-degenerate. We 
described all blowups of non-degenerate $cAx/4$-singularities which
can have non-rational exceptional divisors with small discrepancies 
in section~\ref{S:cAx4}. But in this case all of them are rational.
It follows that only the blowup $\nu$ of $(X,o)$ can have a 
non-rational exceptional divisor. It has the form 
$$E=\{u^2+x^3+\lambda_0xy^4+\mu_0y^6=0\}\subset
\mathbb{P}(2,1,4,3)\,.$$
This is a cone over a curve of genus $g\leqslant 1$. Discrepancy
$a(E,X)$ is equal to $(1/3)(2+1+4+3)-1-2=1/3$. We have proved
\begin{proposition}\label{P:cD32}
There is not more than $1$ non-rational divisor $E$ with discrepancy
$a\leqslant 1$ over the singularity $(X,o)$ or type $cD/3-2$. 
If $X$ is defined by equation~\eqref{E:cD32}, then the non-rational
divisor $E$ is birational to the exceptional divisor of the
blowup $\frac{1}{3}(2,1,4,3)$. It is a cone over a curve of
genus $1$.  
\end{proposition}

\subsubsection{$cD/3-3$}\label{S:cD33}
Consider the singularity $(X,o)$ of type $cD/3-3$, i.~e.,
$$X\simeq\{\varphi=u^2+x^3+y^3+xyz^3\alpha(z^3)+xz^4\beta(z^3)+
yz^5\gamma(z^3)+z^6\delta(z^3)\}\,,$$
where $\alpha(z^3)$, $\beta(z^3)$, $\gamma(z^3)$, 
$\delta(z^3)\in\mathbb{C}\{z^3\}$. Here we additionally assume that
the defining series $\varphi$ is non-degenerate. If $E$ is a 
non-rational divisor with $a(E,X)\leqslant 1$ and $\centr_X(E)=o$,
then, as in section~\ref{S:cAx4}, we can consider $E$ as an 
exceptional divisor of some weighted blowup or pseudo blowup. 
Let $w$ be weight of this blowup. The Newton diagram $\Gamma(f)$ 
is spanned by the monomials $u^2$, $x^3$, $y^3$, $xyz^{3b_1}$,
$xz^{4+3b_2}$, $yz^{5+3b_3}$, $z^{6+3b_4}$, where $b_i\in
\mathbb{Z}_{\geqslant 0}$. Thus if $E$ is non-rational, its
equation $\varphi_w$ contains the monomials $u^2$ and $x^3$, $u^2$ 
and $y^3$, or $x^2$ and $y^3$. Using the condition $a(E)\leqslant 1$,
we come to the following problem: find all primitive vectors
$w\in\mathbb{Z}^4+\frac{1}{3}(1,2,2,0)\mathbb{Z}$ such that \\
either (i) $w_1+w_2+w_3+w_4-1-2w_4\leqslant 1$, $2w_4=3w_1$,
$3w_2\geqslant 2w_4$; \\
or (ii) $w_1+w_2+w_3+w_4-1-2w_4\leqslant 1$, $2w_4<3w_1$,
$3w_2=2w_4$; \\
or (iii) $w_1+w_2+w_3+w_4-1-2w_1\leqslant 1$, $w_1=w_2$, $2w_4>w_1$.
\begin{proposition}\label{P:cD33}
Let the primitive vector $w\in\mathbb{Z}^4+\frac{1}{3}(1,2,2,0)
\mathbb{Z}$ satisfy one of conditions (i), (ii), or (iii).
Then $w$ is one of the following: \\
1) $\frac{1}{3}(5,4,1,6)$; \\
2) $\frac{1}{3}(2,4,1,3)$; \\
3) $\frac{1}{3}(4,5,2,6)$; \\
4) $(2,2,1,3)$.
\end{proposition}
\begin{proof} It is an easy arithmetic calculation.
\end{proof}

Note that weight 4) corresponds to a pseudo blowup, other weights
correspond to weighted blowups. 

The exceptional divisor $E_1$ of the blowup $\nu_1=
\frac{1}{3}(5,4,1,6)$ is defined in $\mathbb{P}(5,4,1,6)$ by the
equation
$$u^2+y^3+\gamma_1yt^8+\delta_2z^{12}=0\,.$$
(Recall that we assume that $E_1$ is non-rational. It follows
that $\alpha_0=\beta_0=$\linebreak
$\beta_1=\gamma_0=\delta_0=\delta_1=0$).
Discrepancy $a(E)=1/3$. The divisor $E$ is a cone over a curve of
genus $1$. 

The exceptional divisor $E_2$ of the blowup $\nu_2=
\frac{1}{3}(2,4,1,3)$ is given in $\mathbb{P}(2,4,1,3)$ by the
equation
$$u^2+x^3+\beta_0xz^4+\delta_0z^6=0\,.$$
Discrepancy $a(E_2)=1/3$ and $E_2$ is again a cone over a curve of
genus $1$.

The exceptional divisor $E_3$ of the blowup $\nu_3=
\frac{1}{3}(4,5,2,6)$ is defined in $\mathbb{P}(4,5,2,6)$ by the
equation 
$$u^2+x^3+\beta_0xz^4+\delta_0z^6=0\,.$$
It follows that $E_3\simeq\{u^2+x^3+\beta_0xz^4+\delta_0z^6=0\}
\subset\mathbb{P}(2,5,1,3)$. It is again a cone over a curve of
genus $1$. Discrepancy $a(E_3)=2/3$.

The exceptional divisor $E_4$ of the blowup $\nu_4=(2,2,1,3)$
is defined as
$$E_4\simeq\{u^2+x^3+y^3+\beta_0xz^4+\delta_0z^6=0\}\subset
\mathbb{P}(2,2,1,3)/G\,,$$
where $G$ is a cyclic group. But the surface $\{u^2+x^3+y^3+
\beta_0xz^4+\delta_0z^6=0\}\subset\mathbb{P}(2,2,1,3)$ has only
rational singularities. According to Lemma~\ref{L:rationality},
the surface $E_4$ is rational.

It is clear that the blowups $\nu_1$ and $\nu_2$, $\nu_1$ and 
$\nu_3$ can not simultaneously be non-rational. But if one of
the blowups $\nu_2$, $\nu_3$ is non-rational, then the other is too. 
\begin{example}
Consider the singularity
$$\{u^2+x^3+y^3+z^6=0\}\subset\frac{1}{3}(1,2,2,0)$$
of type $cD/3-3$ and its blowups $\nu_2$ and $\nu_3$. Their
exceptional divisors  
$$E_2=\{u^2+x^3+z^6=0\}\subset\mathbb{P}(2,4,1,3)$$
$$\text{and }E_3=\{u^2+x^3+z^6=0\}\subset\mathbb{P}(2,5,1,3)$$
are cones over elliptic curves. It is interesting that they are
given by the same equations. But the blowups $\nu_2$ and $\nu_3$ are
not isomorphic since their discrepancies are different: 
$a(E_2)=1/3$, $a(E_3)=2/3$.
\end{example}

Thus there are not more than $2$ non-rational divisors with
discrepancy $a\leqslant 1$ over a singularity of type $cD/3-3$. 

\subsection{Terminal singularities of type $cD/2$}\label{S:cD2}
\subsubsection{$cD/2-1$}\label{S:cD21$}
Consider the singularity $(X,o)$ of type $cD/2-1$, i.~e.,
$$X\simeq\{\varphi=u^2+xyz+x^{2a}+y^{2b}+z^c=0\}\subset
\frac{1}{2}(1,1,0,1)\,,$$
where $a$, $b\geqslant 2$, $c\geqslant 3$. This singularity is
non-degenerate. Thus all divisors with discrepancy $a\leqslant 1$ 
correspond to faces of the Newton diagram $\Gamma(\varphi)$. 
But it is easy to show that all faces produce rational divisors,
hence there are no non-rational divisors with discrepancy
$a\leqslant 1$ over a singularity of type $cD/2-1$.

\subsubsection{$cD/2-2$}\label{S:cD22}
Consider the singularity $(X,o)$ of type $cD/2-2$, i.~e.,
$$X\simeq\{\varphi=u^2+y^2z+\lambda yx^{2a+1}+g(x,z)=0\}\subset
\frac{1}{2}(1,1,0,1)\,.$$
where $\lambda\in\mathbb{C}$, $a\geqslant 1$, $g(x,z)\in
(x^4,x^2z^2,z^3)\mathbb{C}\{x,z\}$. Here we assume that the
series $\varphi$ is non-degenerate. Since the singularity $(X,o)$ is
isolated and the function $g$ is $\mathbb{Z}_2$-invariant, we see
that $g$ contains a monomial of the form $z^{n-1}$. If $n$ is the
minimal integer with this property, then we say that $(X,o)$ is
of type $cD_n/2-2$.

Divisors with discrepancy $a\leqslant 1$ over $(X,o)$ correspond
to faces of the Newton diagram $\Gamma(\varphi)$. In the same way as
in sections~\ref{S:cAx4} and \ref{S:cD33}, we come to the
following problem: find all primitive vectors $w\in
\mathbb{Z}^4+\frac{1}{2}\mathbb{Z}$ such that \\
either (i) $w_1+w_2+w_3+w_4-1-2w_4\leqslant 1$, $2w_2+w_3\geqslant 
2w_4$, $(n-1)w_3\geqslant 2w_4$; \\
or (ii) $w_1+w_2+w_3+w_4-1-2w_2-w_3\leqslant 1$, $2w_4>2w_2+w_3$,
$(n-1)w_3\geqslant 2w_4$. \\
The answer is given by the following
\begin{proposition}\label{P:cD22}
Let the primitive vector $w\in\mathbb{Z}^4+\frac{1}{2}(1,1,0,1)
\mathbb{Z}$ satisfy conditions (i) or (ii).
Then $w$ is one of the following: \\
1) $\frac{1}{2}(1,m,2,m)$, $m=2k-1$, $m\leqslant n-1$; \\
2) $\frac{1}{2}(1,m-2,4,m)$, $m=2k-1$, $m\leqslant 2(n-1)$; \\
3) $\frac{1}{2}(1,m-1,2,m+1)$, $m=2k$, $m\leqslant n-1$; \\
4) $(1,k,2,k)$, $k\leqslant (n-1)/2$; \\
5) $(1,k-1,2,k)$, $k\leqslant n-1$; \\
6) $(1,k-1,1,k)$, $k\leqslant n/2$.
\end{proposition}
\begin{proof} It is an easy arithmetic calculation.
\end{proof}

Blowups $\nu_1=\frac{1}{2}(1,m,2,m)$, $\nu_2=\frac{1}{2}(1,m-2,4,m)$,
$\nu_3=\frac{1}{2}(1,m-1,2,m+1)$ are weighted, and
$\nu_4=(1,k,2,k)$, $\nu_5=(1,k-1,2,k)$, and $\nu_6=(1,k-1,1,k)$ are
pseudo blowups. Actually, only the blowups $\nu_1$, $\nu_3$ (with
discrepancy $a=1/2$), $\nu_4$, and $\nu_6$ (with discrepancy $a=1$) 
can be non-rational.

\begin{example}
Consider the singularity
$$\{u^2+y^2z+z^{2k}+x^{2k}=0\}\subset\frac{1}{2}(1,1,0,1)$$
of type $cD_{2k+1}/2$ and its pseudo blowup $\nu_4=(1,k,1,k)$. 
Assume that the number $k$ is even. Then the affine chart
$U_1=X(\sigma_1,N')$ of the blown up variety 
${\widetilde{\mathbb{C}^4}}_{(1,k,1,k)}$
(for the notation see section~\ref{S:prelim}) is isomorphic to
$$\mathbb{C}^4/\mathbb{Z}_2(1,1-k,-1,1-k)=\frac{1}{2}(1,1,1,1)\,.$$
$$\Tilde X\cap U_1=\{y_{4}^{2}+y_1y_{2}^{2}y_3+y_{3}^{2k}+1=0\}
\subset\frac{1}{2}(1,1,1,1)\,,$$
the exceptional divisor ($y_1=0$)
$$E\simeq\{y_{4}^{2}+y_{3}^{2k}+1=0\}\subset\frac{1}{2}(1,1,1)$$
is a cone over the curve $\{y_{4}^{2}+y_{3}^{2k}+1=0\}\subset
\frac{1}{2}(1,1)$.
It is a hyperelliptic curve of genus $k/2$.
\end{example}

For all the blowups $\nu_1$, $\nu_3$, $\nu_4$, $\nu_6$, the
exceptional divisor is a cone over a hyperelliptic curve
of genus $g\leqslant k-1$ for $\nu_1$; $g\leqslant k$ for $\nu_3$; 
$g\leqslant k/2$ for even $k$ and $g\leqslant (k-1)/2$ for
odd $k$ for $\nu_4$; 
$g\leqslant (k-1)/2$ for odd $k$ and $g\leqslant (k-2)/2$ for even 
$k$ for $\nu_6$. In the latter case the exceptional divisor $E_6$ 
splits onto $2$ components, one of them is rational.

The pairs of blowups $\nu_1$ and $\nu_3$, $\nu_4$ and $\nu_6$ can
not simultaneously be non-rational; the others can. 
\begin{example}
Consider the singularity
$$\{u^2+y^2z+z^{12}+z^6x^6+x^{18}=0\}\subset\frac{1}{2}(1,1,0,1)$$
of type $cD_{13}/2-2$, its weighted blowup $\nu_1=\frac{1}{2}
(1,9,2,9)$, and its pseudo blowup $\nu_4=(1,6,1,6)$.

The exceptional divisors are $E_1$ and $E_4$.
$$E_1=\{u^2+z^6x^6+x^{18}=0\}\subset\mathbb{P}(1,9,2,9)$$
is a cone over a singular curve of genus $2$;
$$E_4=\{u^2+z^{12}+z^6x^6=0\}\subset\mathbb{P}(1,6,1,6)/
\mathbb{Z}_2$$
is a cone over a singular curve of genus $1$.
\end{example}

So, there are not more than $2$ non-rational divisors with
discrepancy $a\leqslant 1$ over a non-degenerate singularity of type 
$cD/2-2$.

\subsection{Terminal singularities of type $cE/2$}\label{S:cE2}
Consider the singularity $(X,o)$ of type $cE/2$, i.~e.,
$$X\simeq\{\varphi=u^2+x^3+g(y,z)x+h(y,z)=0\}\subset
\frac{1}{2}(0,1,1,1)\,,$$
where $g(y,z)\in(y,z)^4\mathbb{C}\{y,z\}$, $h(y,z)\in
(y,z)^4\mathbb{C}\{y,z\}\setminus(y,z)^5\mathbb{C}\{y,z\}$. We
assume that the series $\varphi$ is non-degenerate. In addition,
by permutation of coordinates $y$ and $z$ if necessary, we can
suppose that $y^4$, $y^3z$, or $y^2z^2\in h(y,z)$. The argument here
is similar to that of sections~\ref{S:cAx4}, \ref{S:cD33}, and 
\ref{S:cD22}, so we just formulate the final results.

For $cE/2$-singularities, the non-rational divisors again
can be represented as exceptional divisors of certain weighted
blowups and pseudo blowups. All these divisors are cones. In the
following proposition we list all possible non-rational
blowups, discrepancies of their exceptional divisors, and genuses
of the corresponding curves. 
\begin{proposition}\label{P:cE2}
(cf. \cite{Hay}, \S10) Let $E$ be a non-rational divisor over the 
singularity $(X,o)$ such that $\centr_X(E)=o$ and 
$a(E,X)\leqslant 1$. Then $E$ is birational to the exceptional 
divisor of one of the following blowups. \\
1) $\nu_1=\frac{1}{2}(2,3,1,3)$, $a=1/2$, $g=1$; \\
2) $\nu_2=\frac{1}{2}(2,1,3,3)$, $a=1/2$, $g=1$; \\
3) $\nu_3=\frac{1}{2}(4,3,1,5)$, $a=1/2$, $g=1$; \\
4) $\nu_4=\frac{1}{2}(4,3,1,7)$, $a=1/2$, $g\leqslant 3$; \\
5) $\nu_5=\frac{1}{2}(6,5,1,9)$, $a=1/2$, $g=1$; \\
6) $\nu_6=(2,2,1,3)$, $a=1$, $g=1$; \\
7) $\nu_7=(3,2,1,4)$, $a=1$, $g=1$.
\end{proposition}

Note that curve for the blowup $\nu_4$ is not necessarily 
hyperelliptic.
\begin{example}
$$\{u^2+x^3+y^3+z^{12}=0\}\subset\frac{1}{2}(0,1,1,1)\,.$$
The exceptional divisor of the weighted blowup $\nu_4$ is given by
the equation
$$\{x^3+y^4+z^{12}=0\}\subset\mathbb{P}(4,3,1,7)\,.$$
It is a cone over a non-hyperelliptic curve of genus $3$.
\end{example}

Only the following pairs of blowups can simultaneously be
non-rational: $\nu_1$ and $\nu_2$, $\nu_1$ and $\nu_6$.
\begin{example}
$$\{u^2+x^3+y^2z^2+y^6+z^6=0\}\subset\frac{1}{2}(0,1,1,1)\,.$$
The exceptional divisor of the blowup $\nu_1$ is
$$\{u^2+x^3+z^6=0\}\subset\mathbb{P}(2,3,1,3)\,,$$
the exceptional divisor of the blowup $\nu_2$ is
$$\{u^2+x^3+y^6=0\}\subset\mathbb{P}(2,1,3,3)\,.$$
Both of them are cones over elliptic curves.
\end{example}

So, there are not more than $2$ non-rational divisors with
discrepancy $a\leqslant 1$ over a non-degenerate singularity
of type $cE/2$.

\end{document}